\documentclass{amsart}

\newtheorem{theorem}{Theorem}
\newcommand{\bt}{\begin{theorem}}
\newcommand{\et}{\end{theorem}}

\newtheorem{lemma}{Lemma}
\newcommand{\bl}{\begin{lemma}}
\newcommand{\el}{\end{lemma}}

\newtheorem{corollary}{Corollary}
\newcommand{\bc}{\begin{corollary}}
\newcommand{\ec}{\end{corollary}}

\newtheorem{problem}{Problem}
\newcommand{\bp}{\begin{problem}}
\newcommand{\ep}{\end{problem}}

\newcommand{\Z}{\ensuremath{\mathbb Z}}
\newcommand{\N}{\ensuremath{\mathbb N}}
\newcommand{\A}{\ensuremath{\mathcal A}}
\newcommand{\B}{\ensuremath{\mathcal B}}
\newcommand{\U}{\ensuremath{\mathcal U}}

\begin{document}
\large

\title[Perfect difference sets]{Perfect difference sets constructed from
Sidon sets}

\author{Javier Cilleruelo} \thanks{The work of J.C. was supported by Grant MTM 2005-04730 of MYCIT (Spain)}
\address{Departamento de Matem{\'a}ticas\\
Universidad Aut{\'o}noma de Madrid\\
28049 Madrid, Spain}
\email{franciscojavier.cilleruelo@uam.es}

\author{Melvyn B. Nathanson}
\thanks{The work of M.B.N. was supported in part by grants from the
NSA Mathematical Sciences Program and the PSC-CUNY Research Award
Program.}
\address{Department of Mathematics\\
Lehman College (CUNY)\\
Bronx, New York 10468}
\email{melvyn.nathanson@lehman.cuny.edu}

\begin{abstract}
A set $\A$ of positive integers is a \emph{perfect difference set}
if every nonzero integer  has an unique representation as the
difference  of two elements of $\A$. We construct  dense {\em
perfect difference  sets} from dense Sidon sets.  As a consequence
of this new approach we prove that there exists   a perfect
difference set $\A$  such that
\[
A(x) \gg x^{\sqrt{2}-1-o(1)}.
\]

Also we prove that there exists a \emph{perfect difference set}
$\A$ such that \linebreak $\limsup_{x\to \infty}A(x)/\sqrt x\ge
1/\sqrt 2$.
\end{abstract}

\maketitle

\section{Introduction}
Let \Z\ denote the integers and \N\ the positive integers.
For nonempty sets of integers \A\ and \B, we define the \emph{difference set}
\[
\A - \B = \{a-b : a\in \A \text{ and } b\in \B \}.
\]
For every integer $u$, we denote by
$d_{\A,\B}(u)$  the number of pairs $(a,b)\in \A \times \B$ such that
$u=a-b.$
Let $d_{\A}(u)$  the number of pairs $(a,a')\in \A \times \A$ such that
$u=a-a'.$
The set $\A$ is a \emph{perfect difference set} if
$d_{\A}(u) = 1$ for every integer $ u \neq 0$.   Note that  \A\ is a perfect difference set if and only $d_{\A}(u) = 1$ for every positive integer $u$.
For perfect difference sets, a simple counting argument shows that $$A(x)\ll x^{1/2},$$ where the
{\em counting function} $A(x)$ counts the number of positive
elements of $\A$ not exceeding $x$.

It is not completely obvious that perfect difference sets exist, but
 the greedy algorithm  produces \cite{lev04}  a perfect difference
 set $\A\subseteq \N$  such that $$A(x)\gg x^{1/3}.$$
At the Workshop on Combinatorial and Additive Number Theory (CANT
2004) in New York in May, 2004, Seva Lev (see also
\cite{lev04}) asked if there exists a perfect difference set
$\A$ such that \[A(x) \gg x^{\delta} \text{ for
some } \delta
> 1/3.\]
We answer this questions affirmatively  by  constructing perfect difference
sets from classical Sidon sets.

 We say that a set $\B$ is a Sidon set if
 $d_{\B}(u)\le 1$ for all integer $u\ne 0$.

\bt  \label{pds}
For every Sidon set $\mathcal B$ and every function $\omega(x)\to \infty$, there exists a perfect difference set $\A\subseteq \N$ satisfying
 \[A(x)\ge  B(x/3)-
 \omega(x).\]
\et

It is a difficult problem to construct dense infinite Sidon sets.
 Ruzsa~\cite{ruzs98}  proved that there exists a Sidon set $
\mathcal B$ with $  B(x) \gg x^{\sqrt{2}-1-o(1)}.$ The following
result follows easily.

\bt
There exists a perfect
difference set $\A\subseteq \N$  such that
\[A(x)\gg x^{\sqrt 2-1+o(1)}.\]
\et

Erd\H os~\cite{stoh55} proved that the lower bound $A(x)\gg x^{1/2}$ does not hold for any Sidon set $\A$, and so does not hold for perfect difference sets.  However, Kr{\"u}ckeberg~\cite{kruc61} proved that there exists a  Sidon set $\B$
such that $$\limsup_{x\to \infty}\frac{B(x)}{\sqrt{x}}\ge \frac
1{\sqrt 2}.$$   We extend this result to perfect difference sets.

\bt\label{kpds}There exists a perfect difference set $\A\subset \N$ such that
\begin{equation*}
\limsup_{x\to \infty}\frac{A(x)}{\sqrt x}\ge \frac 1{\sqrt 2}.
\end{equation*}
\et

Notice that an immediate application of Theorem~\ref{pds} to
Kr{\"u}ckeberg's result would give only
$\limsup_{x\rightarrow\infty}A(x) x^{-1/2} \geq 1/\sqrt 6.$

\section{Proof of Theorem~\ref{pds}}
\subsection{Sketch of the proof}
The strategy of the proof is the following:
\begin{enumerate}
    \item[$\bullet$] Modify any dense Sidon set $\mathcal B$ given by dilating it by $3$ and  removing  a
suitable {\em thin} subset of $3\ast \B$.
    \item[$\bullet$] Complete  the remainder set
    $\B_0=(3\ast \B)\setminus \{ \text{removed set} \}$ with a subset of the  elements of a very sparse
    sequence $\U=\{ u_s\}$ by adding, if $k$ has not appeared
    yet in the difference set, two elements $u_{2k},u_{2k+1}$ in the $k$-th step
    such that $u_{2k+1}-u_{2k}=k$.
   \end{enumerate}

\subsection{The auxiliary sequence $\U$}

For any strictly increasing function $g:\N \to \N$ and  $k\ge 1$, we define integers $u_{2k}$ and $u_{2k+1}$ by
\[
\begin{cases}u_{2k}&= 4^{g(k)}+\epsilon_k \\
  u_{2k+1}&= 4^{g(k)}+\epsilon_k+k \end{cases}
  \]
where
\[
\epsilon_k=\begin{cases} 1\ \text{ if }
k\equiv 2\pmod 3
\\ 0\ \text{ otherwise}\end{cases}.
\]
For all positive integers $k$ we have
\[
u_{2k+1}-u_{2k} = k.
\]
Let  $ \U_k=\{ u_{2k},u_{2k+1}\}$ and
$\U_{< \ell}=\bigcup_{k<\ell}\U_k$.
 It will be useful to state some properties
of the sequence $\U=\{ u_i\}_{i=2}^{\infty}$.

\bl \label{propU1} The sequence $\U=\{u_i\}_{i=2}^{\infty}$ satisfies the following
properties:\begin{enumerate}
    \item[(i)] For all $i \geq 2$, $u_i\not \equiv 0\pmod 3$.
    \item[(ii)] For all $k \geq 2$, for $u\in \U_k$, and for
    all $u',u'',u'''\in \U_{<k}$, we have $u+u' > u''+u'''$.
    \item[(iii)] If $k \geq 2$, $u\in \U_k$, and
    $u'\in \U_{<k}$, then $u-u'>u/2$.
\end{enumerate}
\el

\begin{proof}
(i) By construction.

(ii) Since $g(k)$ is strictly increasing we have $k \leq g(k)$
and so
\[
4k < 4^k \leq 4^{g(k)}
\]
for all $k \geq 2$.  It follows that
\begin{align*}
u''+u''' & \le 2u_{2k-1} \le 2(4^{g(k-1)}+k) \\
& \le 2\left(4^{g(k)-1}+k\right)  \le \frac{4^{g(k)}}{2}+2k \\
& \leq 4^{g(k)} \leq u < u+u'.
\end{align*}
%

(iii) For $k \geq 2$ we have
\begin{align*}
u'  & \leq 4^{g(k-1)}+ (k-1) +\epsilon_{k-1}\\
&  \leq 4^{g(k)-1}+k \\
& <  2\cdot 4^{g(k)-1} = \frac{4^{g(k)}}{2} \\
& \leq u/2
\end{align*}
and so $u-u'>u/2$.
\end{proof}

\subsection{ Construction of the Sidon set $\B_0$}

Take a Sidon set $\mathcal B$ and consider the set
$\B' =3\ast \mathcal B=\{ 3b : b\in \mathcal B\}$.
Then $\B'$ is a Sidon set such that $b\equiv 0 \pmod 3$ for all  $b\in \B'$
and $B'(x)= B\left(\frac{x}{3}\right )$.

The set $\B_0$ will be the set $\B'=3\ast \B$ after we remove all the
elements $b\in \B'$ that satisfy at least one of the followings conditions:
\begin{description}
    \item[c1] $b=u-u'+b'$ for some $b'\in \B',\ b>b'$ and
    $u,u'\in \U $ such that
    $u\in \U_r,\ u'\in \U_{<r}$
    for some $r$.
    \item[c2] $b=u+u'-b'$ for some $b'\in \B,\ b\ge b'$ and
    $u,u'\in \U$.
    \item[c3] $b=u+u'-u''$ for some $u\in \U_r$, $u' \in \U$, and $u''\in \U_{<r}$ with $u'\le u $.
    \item[c4] $|b-u_i|\le i$ for some $u_i \in \U$.
\end{description}

\subsection{ The inductive step}
We shall construct the set $\A$ in Theorem~\ref{pds} by adjoining
terms to the \emph{nice} Sidon set $\B_0$  obtained above. More
precisely, the sequence $\A$ satisfying the conditions of the
theorem will be
$$\A=\bigcup_{k=0}^{\infty}\A_k $$ where $\A_k$ will be defined by
$\A_{0}=\B_0$ and for, $k\ge 1$,
\[
\A_{k}=
\begin{cases}
\A_{k-1}\cup \U_k   &\text{ if } k\not \in \A_{k-1}-\A_{k-1}  \\
\A_{k-1} & \text{ otherwise.}
\end{cases}
\]

\bl
For every positive integer $k$ we have
\[
[-k,k]\subseteq \A_k-\A_k
\]
and so
\[
d_{\A}(n) \geq 1
\]
for all integers $n$.
\el
\begin{proof}
Clear.
\end{proof}

\subsection{ $\A$ is a Sidon set}
First we state two lemmas.

\bl\label{union}
Let $A_1$ and $A_2$ be nonempty disjoint sets of integers and let $A=A_1\cup A_2$  For every integer $n$ we have
\[
d_A(n)=d_{A_1}(n)+d_{A_2}(n)+d_{A_1,A_2}(n)+d_{A_2,A_1}(n),
\]
where $$d_{A_i,A_j}(n)=\# \{ (a,a')\in A_i\times A_j,\ a-a'=n\}.$$
\el
\begin{proof}
This follows from the identity
\[
(A_1\cup A_2)\times (A_1\cup A_2)=(A_1\times A_1) \cup (A_2\times
A_2) \cup (A_1\times A_2)\cup (A_2\times A_1).
\]
\end{proof}

\bl\label{A-U} If  $n\in \A_{r-1}-\U_r$ then
\begin{enumerate}
    \item [(i)] $|n| >r$, and so
    $d_{\U_r,\A_{r-1}}(r)=d_{\A_{r-1},\U_r}(r)=0$.
    \item [(ii)] $d_{\A_{r-1}}(n)=0$.
    \item [(iii)] $d_{\U_r,\A_{r-1}}(n)=0$.
\end{enumerate}
\el

\begin{proof} Write $n=a-u,$ where $ a\in \A_{r-1}$ and $u\in \U_r = \{u_{2r}, u_{2r+1}\}$.

(i) If $a=b\in \B_0$ we have that $|b-u| > 2r > r$ because, by condition~(c4),  we have removed all elements $b$ from $\B$ such that
$|b-u_i|\le i$.

If $a=u'\in \U_{<r}$ then we apply Lemma~\ref{propU1} (iii) to
conclude that
\[
|u'-u|>\frac{u}{2} \ge\frac{4^{g(r)}}2>r.
\]

(ii) Since $\A_{r-1}\subseteq \B_0\cup \U_{<r}$, it follows that
\[
d_{\A_{r-1}}(n)\le d_{\B_0\cup\U_{<r}}(n)\le
d_{\B_0}(n)+d_{\U_{<r}}(n)+d_{\B_0,\U_{<r}}(n)+d_{\U_{<r},
\B_0}(n).
\]

If $a=b\in \B_0$, then $n = b-u$ and
\begin{enumerate}
\item
$b \equiv 0\pmod{3}$ but $u \not\equiv 0 \pmod{3}$,  hence $b-u \not\equiv 0 \pmod{3}$ and $d_{\B_0}(b-u)=0$ (by Lemma \ref{propU1}~(i)),
\item
$d_{\U_{<r}}(b-u)=0$ (by condition~(c3)),
\item
$d_{\B_0,\U_{<r}}(b-u)=0$ (by condition (c1)),
\item
$d_{\U_{<r}, \B_0}(b-u)=0$ (by condition~(c2)).
\end{enumerate}

If $a=u'\in \U_{<r}$, then $n = u' - u$ and
\begin{enumerate}
\item $d_{ \B_0}(u'-u)=0$ (by condition (c1)), \item
$d_{\U_{<r}}(u'-u)=0$ (by  Lemma \ref{propU1} (ii)), \item If
$u'-u=b-u''$ with $u''\in \U_{<r}$, then Lemma~\ref{propU1} (iii)
implies that $0<b=u'+u''-u \leq 0$, and so
$d_{B_0,\U_{<r}}(u'-u)=0$. \item $d_{\U_{<r},B_0}(u'-u)=0$ (by
condition (c3)).
\end{enumerate}

(iii) Again, since $\A_{r-1}\subseteq \B_0\cup \U_{<r}$ we have
that
\[
d_{\U_r,\A_{r-1}}(n)\le d_{\U_r,\B_0}(n)+d_{\U_r,\U_{<r}}(n).
\]
If $a=b\in \B_0$ then $d_{\U_r,\B_0}(b-u)=0$ (by condition (c2))
and $d_{\U_r,\U_{<r}}(b-u)=0$ (by condition (c3)).

If $a=u'\in \U_{<r}$ then $d_{\U_r,\B_0}(u'-u)=0$ (by condition
(c3)).  Finally, we have $d_{\U_r,\U_{<r}}(u'-u)=0$, since if
$u'-u=u''-u''',\ u''\in \U_r,\ u'''\in \U_{<r}$, then
$0>u'-u=u''-u'''>0$. This completes the proof.
\end{proof}

\bl For every positive integer $n$ we have
$$d_{\A}(n)\le 1$$
and so \A\ is a perfect difference set.
\el

\begin{proof}
We will use induction to prove that, for every $r \geq 0$,
\[
d_{\A_r}(n)\le 1\qquad \text{ for every nonzero integer } n.
\]
This is true for $r=0$ because $\A_{0}=\B_0$ is a subset of a
Sidon set.

We assume that the statement is true for $r-1$ and shall prove it
for $r$.

If $d_{\A_{r-1}}(r)=1$ then $\A_r=\A_{r-1}$ and there is nothing
to prove.  Suppose that  $d_{\A_{r-1}}(r)=0$, and so
$\A_r=\A_{r-1}\cup \U_r$.  Since we have added two new elements
$u_{2r}, u_{2r+1}$ to $\A_{r-1}$, it is possible that there are
{\it new} representations of a positive integer $n$ so that
$d_{A_r}(n)>1$. We shall prove that this cannot happen.

By Lemma \ref{union}, we can write
\[
d_{\A_r}(n)=d_{\A_{r-1}}(n)+d_{\U_r}(n)+d_{\A_{r-1},\U_r}(n)+d_{\U_r,\A_{r-1}}(n)
\]
If $n=r$, then Lemma \ref{A-U} (i) and the relation
$u_{2r+1}-u_{2r}=r$ imply that
\[
d_{\A_r}(r)=d_{\A_{r-1}}(r)+d_{\U_r}(r)+d_{\A_{r-1},\U_r}(r)+d_{\U_r,\A_{r-1}}(r)=
0+1+0+0=1
\]
If $n\ne r$, then
\[
d_{\A_r}(n)=d_{\A_{r-1}}(n)+d_{\A_{r-1},\U_r}(n)+d_{\U_r,\A_{r-1}}(n).
\]

If $n\in \A_{r-1}-\U_r$ (the case $n\in \U_r-\A_{r-1}$ is
similar), then we can write
$$
n=a-u \text{ where } a\in \A_{r-1},\ u\in \U_r.
$$
Applying Lemma \ref{A-U} (ii) and Lemma \ref{A-U} (iii), we obtain
$$
d_{\A_r}(n)=d_{\A_{r-1},\U_r}(n).
$$
If $d_{\A_{r-1},\U_r}(n)\ge 2$, then there exist $a,a'\in
\A_{r-1}$ such that $a-u_{2r}=a'-u_{2r+1}$.  This implies that
\[
a'-a=u_{2r+1}-u_{2r}=r \in \A_{r-1}-\A_{r-1}
\]
which is false, so $d_{\A_r}(n)=d_{\A_{r-1},\U_r}(n)\le 1$.

If $n\not \in \left ( \A_{r-1}-\U_r\right )\cup \left (
\U_r-\A_{r-1}\right )$ then
\[
d_{\A_r}(n)=d_{\A_{r-1}}(n)\le 1.
\]
This completes the proof.
\end{proof}

\subsection{The counting function $A(x)$}
We have
\[
A(x)\ge
B_0(x)=B(x)-R(x)=B(x/3)-R(x)
\]
where $R=R_1\cup
R_2\cup R_3\cup R_4$ and $R_i$ denotes the set of elements of $B$
removed by condition $(c_i)$, $i=1,2,3,4$.

\bl
Let $U(x )$ denote the counting function of the set $\U$.
For the sets $R_1,R_2,R_3,R_4$ defined above, we have
\begin{enumerate}
    \item[(i)] $R_1(x)\le U^2(2x)$.
    \item[(ii)] $R_2(x)\le U^2(2x)$.
    \item[(iii)] $R_3(x)\le U^3(2x)$
    \item[(iv)] $R_4(x)\le 2U^2(2x) + U(2x)$.
\end{enumerate}
\el

\begin{proof}
(i) We have
\[
R_1(x)=\# \{  b\in B : b\le x  \text{ and $b$ satisfies condition (c1)} \}.
\]
Because $B$ is a Sidon set,  for every pair of integers  $u,u' \in
\U $ there exists at most one pair of integers $b,b'\in \B$ such
that $b-b'=u-u'$.  The condition $x\ge b>b'$ implies that
$0<u-u'\le x$.  On the other hand Lemma \ref{propU1} (iii) implies
that $u-u'>u/2$ and so $u<2x$ and
\[
R_1(x)\le \# \{ (u,u'),\ u'<u, \ u< 2x\}\le U^2(2x).
\]

(ii) Again, because $B$ is a Sidon set, for every pair $u,u' \in \U$
there exists at most one pair $b,b'\in \B$ such that $b+b'=u+u'$. The
condition $x\ge b \geq b'$ implies $u,u'\le 2x$ and so
\[
R_2(x)\le \# \{ (u,u') \in \U \times \U : u\le 2x, u'\le 2x\} \le
U^2(2x).
\]

(iii) If  $u\in \U_r,\ u''\in \U_{<r}$, then Lemma \ref{propU1}
(iii) implies that $b=u+u'-u''>u-u'' > u/2$ and so
\begin{align*}
R_3(x) & = \# \{  b\in B : b\le x  \text{ and $b$ satisfies condition (c3)} \} \\
& \le \# \{ (u,u',u'') \in \U \times \U \times \U :  u<2x, u''<u, u'\le u \}  \\
& \le U(2x)^3.
\end{align*}

(iv) We have
\begin{align*}
R_4(x)  & = \# \{  b\in B : b\le x \text{ and } |b-u_i| \leq i \text{ for some $u_i \in \U$} \}\\
& \le \# \{ n \in \N : n\le x \text{ and }  |n-u_i|\le i \text{ for some } i\}.
\end{align*}
If $n \leq x$ and $|n-u_i|\le i$, then $u_i\le n+i\le x+i$.
Since $u_2 = 4^{g(1)} \geq 4$,  $u_3 = 4^{g(1)+1} \geq 16$, and, for $i \geq 4$,
\[
u_i \geq 4^{g((i-1)/2)} \geq 4^{(i-1)/2} = 2^{i-1} \geq 2i.
\]
Therefore, $u_i \leq x+i \leq x+u_i/2$ and so $u_i \leq 2x$.
It follows that $i\le U(2x)$ and so
\begin{align*}
R_4(x) & \le \# \{ n\le x : |n-u_i|\le U(2x) \text{ and } u_i\le 2x\} \\
& \leq (2U(2x)+1)U(2x) = 2U(2x)^2 +U(2x).
\end{align*}
This completes the proof of the lemma.
\end{proof}

Finally, given any function $\omega(x)\to \infty$ we have that
$$
A(x)\ge B(x/3)-\left (U(2x)^3+4U^2(2x)+U(2x)\right )\ge
B(x/3)-\omega(x)
$$
for any function $g:\N \rightarrow \N$ and sequence $\U$ growing fast enough.
This completes the proof of Theorem~\ref{pds}.

\section{Proof of Theorem~\ref{kpds}}
\bl \label{sidon} If $C_1$ and $C_2$ are Sidon sets such that
$(C_i-C_i)\cap (C_j-C_j)=\{0\}$, $(C_i+C_i)\cap
(C_j+C_j)=\emptyset$ and $(C_i+C_i-C_i)\cap C_j=\emptyset$ for \
$i\ne j$, then $C_1\cup C_2$ is a Sidon set.
\el
\begin{proof}
Obvious.
\end{proof}

\bl For each odd prime $p$ there exist a Sidon set $\B_p$ such
that
\begin{enumerate}
    \item[(i)] $\B_p\subseteq [1,p^2].$
    \item[(ii)] $ (\B_p-\B_p)\cap [-\sqrt p,\sqrt p]=\emptyset$.
    \item[(iii)] $|\B_p|>p-2\sqrt p.$
\end{enumerate}
\el

\begin{proof}
Ruzsa~\cite{ruzs93} constructed, for each prime $p$, a Sidon set $R_p\subseteq [1,p^2-p]$ with $|R_p|=p-1$.   We consider the subset $\B_p$ of $R_p$ that we obtain by removing all elements $b\in R_p$ such that $0<|b-b'|\le \sqrt p$ for some $b'\in R_p$.  Since $R_p$ is a
Sidon set, it follows that we have removed at most $\sqrt{p}$ elements from $R_p$, and so  $|\B_p|\ge |R_p|-\sqrt p=p-\sqrt p-1>p-2\sqrt p$.
\end{proof}

\begin{proof}[Proof of Theorem~\ref{kpds}]
We shall construct an increasing sequence of finite set $A_1 \subseteq A_2 \subseteq A_3 \subseteq \cdots $ such that $\A = \cup_{k=1}^{\infty} A_k$
is a perfect difference set satisfying Theorem~\ref{kpds}.

In the following, $l_k$ will denote the largest integer in the set
$A_{k-1}$, and $p_{k}$ the least prime greater than $4l_k^2$. Let
\[
A_1=\{0,1\}.
\]
Then $l_2=1$ and $p_2=5.$  We define
\begin{equation*}
A_k=\begin{cases} A_{k-1}\cup \left (\B_{p_k}+p_k^2+2l_k\right ) &
\text{ if } k\in A_{k-1}-A_{k-1}\\
A_{k-1}\cup \left (\B_{p_k}+p_k^2+2l_k\right )\cup \{
4p_k^2,4p_k^2+k\} &  \text{ otherwise}.
\end{cases}
\end{equation*}
We shall prove that the set $\A=\cup_{k=1}^{\infty}A_k$ satisfies the theorem.

By construction, $[1,k]\subseteq A_k-A_k$ for every positive integer $k$ and so $\A-\A = \Z$.

We must prove that $A_k$ is a Sidon set for every $k\ge 1$.

This is clear for $k=1$. Suppose that $A_{k-1}$ is a Sidon set.
Let  $C_1=A_{k-1}$ and $C_2=\B_{p_k}+p_k^2+2l_k$. We shall show
that
\[
C_1 \cup C_2 = A_{k-1}\cup (\B_{p_k}+p_k^2+2l_k)
\]
is a Sidon set. Notice that
\[
C_1-C_1\subseteq [-l_k,l_k]\subseteq [-\sqrt{p_k},\sqrt{p_k}]
\]
\[
C_2-C_2=\B_{p_k}-\B_{p_k}
\]
\[
 [ -\sqrt{p_k},\sqrt{p_k}] \cap
(\B_{p_k}-\B_{p_k})=\{ 0\}.
\]
Then $$(C_1-C_1)\cap (C_2-C_2)=\{ 0\}.$$

 Notice also that if
$x\in C_2+C_2$ then $x\ge 2p_k^2+4l_k$, but $C_1+C_1\subset
[1,2l_k]$. Then $$ (C_1+C_1)\cap (C_2+C_2)=\emptyset. $$ If $x\in
(C_1+C_1-C_1)$, then $x\le 2l_k$, but if $x\in C_2$, then
$x>2l_k$. Thus,
$$(C_1+C_1-C_1)\cap C_2=\emptyset.$$
If $x\in C_2+C_2-C_2$, then $x\ge
2(p_k^2+2l_k+1)-(p_k^2+p_k^2+2l_k)=2l_k+1$, and if $x\in C_1$,
then $x\le l_k$. Therefore,
$$(C_2+C_2-C_2)\cap C_1=\emptyset.$$ç
Then $A_{k-1}\cup (\B_{p_k}+p_k^2+2l_k)$ is a Sidon set.

\

Now we must distiguish two cases:

If $k\in A_{k-1}-A_{k-1}$ then $A_k=A_{k-1}\cup
(\B_{p_k}+p_k^2+2l_k)$ and we have proved that it is a Sidon set.

If $k\not \in A_{k-1}-A_{k-1}$ then $A_k=A_{k-1}\cup
(\B_{p_k}+p_k^2+2l_k)\cup \{ 4p_k^2,4p_k^2+k\}$ and we have to
prove that it is also a Sidon set. In this case we take
$C_1=A_{k-1}\cup (\B_{p_k}+p_k^2+2l_k)$ and $C_2=\{
4p_k^2,4p_k^2+k\}$.  We can write
\begin{eqnarray*}
C_1-C_1  = & (A_{k-1}-A_{k-1}) \cup (\B_{p_k}-\B_{p_k})  \\
& \cup  (A_{k-1}-(\B_{p_k}+p_k^2+2l_k)) \\
& \cup \left ((\B_{p_k}+p_k^2+2l_k)-A_{k-1}\right ).
\end{eqnarray*}

If $x\in \left (A_{k-1}-(\B_{p_k}+p_k^2+2l_k)\right )\cup \left
((\B_{p_k}+p_k^2+2l_k)-A_{k-1}\right )$, then $|x|\ge
p_k^2+l_k>k$.

If $x\in (\B_{p_k}-\B_{p_k})$ then $x=0$ or $|x|>\sqrt{p_k}>k$,
then, since $C_2-C_2=\{-k,0,k\}$, we have
$$(C_1-C_1)\cap (C_2-C_2)=\{ 0\}.$$

On the other hand if $x\in C_2+C_2$ then $x\ge 8p_k^2$ but
\[
C_1+C_1\subset [1,2((p_k^2-p_k)+p_k^2+2l_k)]\subset [1,4p_k^2].
\]
Then $$(C_1+C_1)\cap (C_2+C_2)=\emptyset.$$

If $x\in C_1+C_1-C_1$ then $x\le 3p_k^2+2l_k<4p_k^2$. Thus,
$$(C_1+C_1-C_1)\cap C_2=\emptyset.$$ Also we have that
$C_2+C_2-C_2=4p_k^2+\{ -k,0,k,2k\}$, but if $x \in C_1$ we have
that $x<2p_k^2+2l_k<2p_k^2+\sqrt{p_k}<4p_k^2-k$. Thus
$$(C_2+C_2-C_2)\cap C_1=\emptyset.$$

\smallskip

To finish the proof of the theorem note that
\begin{eqnarray*}\limsup_{x\to \infty}\frac{\A(x)}{\sqrt x}\ge \limsup_{k\to
\infty}\frac{\A(2p_k^2-p_k+l_k)}{\sqrt{2p_k^2-p_k+l_k }}\ge
\\\limsup_{k\to
\infty}\frac{|\B_{p_k}|}{\sqrt{2p_k^2-p_k+l_k }}\ge \limsup_{k\to
\infty}\frac{p_k-2\sqrt{p_k}}{ \sqrt{2p_k^2-p_k+\sqrt{p_k}/2 }
}=\frac 1{\sqrt 2}.\end{eqnarray*}
\end{proof}

\section{Remarks and Open problems}
\subsection{The sequence $t(\A)$ associated to a perfect difference
set}

Any translation of a perfect difference
set intersects to itself in exactly one element, and so  we can
define, for every perfect difference set $\A$, a  sequence $t(\A)$
whose elements are given by $t_n=\A \cap (\A-n)$ for all $n\ge 1$.
The sequence $t_n$ is  very irregular, but the greedy algorithm
used in \cite{lev04} generates a perfect difference set such that
$t_n\ll n^3$.  Our method generates a dense Sidon set $\A$, but
gives a very poor upper bound for the sequence $t_n$.

\bp
Does there exists  perfect difference set such that
$t_n=o(n^3)$?
\ep

%
%

\subsection{Sidon sets included in  perfect difference sets}

We have proved that any Sidon set can be perturbed slightly to become a subset of a perfect difference set.   Every subset of a perfect difference set is a Sidon set.  It is natural to ask if \emph{every} Sidon set is a subset of a perfect difference set.  The answer is negative.  To construct a counterexample, we take a perfect
difference set $\A$ and consider the set $\B= 2\ast \A = \{2a: a\in \A\}$. The set $\B$ has the following properties:
\begin{itemize}
    \item[(i)] $\B$ is a Sidon set.
    \item[(ii)] If $n$ is an even integer not in \B, then $\B \cup \{ n\}$ is not a Sidon set.
    \item[(iii)] If $m$ and $m'$ are distinct odd integers not in \B, then $\B \cup \{ m,m'\}$ is not a Sidon set.
\end{itemize}
The Sidon set $\B$ is not a subset of a perfect difference set.
Since this construction is rather artificial, we wonder if almost all
Sidon sets are subsets of perfect difference sets.

\bp
Determine when a Sidon set is a subset of a perfect difference set.
\ep

\subsection{Perfect $h$-sumsets}
Let $\A$ be a set of of integers.  For every integer $u$, we denote by $r^h_{\A}(u)$ the number of
$h$-tuples  $(a_1,\dots ,a_h)\in \A ^h,$ such that $$a_1\le \cdots \le a_h$$
and
$$a_1+\cdots +a_h=u.$$

We say that $\A$ is a {\em perfect $h$-sumset} or a \emph{unique
representation basis of order $h$} if $r^h_{\A}(u)=1$ for every
integer $u$. Nathanson \cite{nath05} proved that for every $h\ge
2$ and for every function $f:\Z \to \N_0 \cup \{\infty\}$ such
that $\limsup_{|u|\to \infty} f(u)\ge 1$ there exists a set of
integers $\A$ such that
$$r_{\A}^h(u)=f(u)$$ for every integer $u$. In particular, the {\em
perfect $h$-sumsets} correspond to the representation function
$f\equiv 1$. Nathanson's construction produces a {\em perfect
$h$-sumset} $\A$ with $$A(x)\gg x^{1/(2h-1)}$$ and he asked for
denser constructions.

It is easy to modify our approach to get a perfect $2$-sumset $\A$
with $A(x)\gg x^{\sqrt2 -1+o(1)}$. But for $h\ge 3$ our method
cannot be adapted easily, and a more complicated construction is
needed. We shall study perfect $h$-sumsets in a forthcoming
paper~\cite{cill-nath06}.

\subsection{Sums and differences}
Let \A\ be a set of integers.  For every integer $u$, we denote by $d_{A}(u)$ and $s_{A}(u)$ the number of solutions of
$$u=a-a' \text{ with } a,a'\in \A
$$
and
$$u=a+a' \text{ with $a,a'\in \A$ and $ a\le a'$,}
$$
respectively. We say that $\A$ is a \emph{perfect difference
sumset} if $d_{\A}(n)=1$ for all  $n\in \N $ and if $s_{\A}(n)=1$
for all $n\in \Z$.

We can extend Theorem~\ref{pds} and Theorem~\ref{kpds} to perfect
difference sumsets.  Then it is a natural to ask if, for any two
functions $f_1:\N\to \N$ and $ f_2:\Z \to \N$, there exists a set
$\A$ such that $d_{\A}(n)=f_1(n)$ for all $n \in \N$ and
$s_{\A}(n)=f_2(n)$ for all $n \in \Z$.  (Note that perfect
difference sumsets correspond to the functions $f_1\equiv 1$ and
$f_2\equiv 1$.) It is not difficult to guess that the answer is
no. For example, if $s_{\A}(n)=2$ for infinitely many integers
$n$, it is easy to see that $d_{\A}(n)\ge 2$ for infinitely many
integers $n$.

\bp Give general conditions for functions $f_1$ and $f_2$ to
assure that there exists a set $\A$ such that $d_{\A}(n)\equiv
f_1(n)$
 and $s_{\A}(n)\equiv f_2(n)$.
 \ep

Is the condition $\liminf_{u\to \infty}f_1(u)\ge 2$ and
$\liminf_{|u|\to \infty}f_2(u)\ge 2$ sufficient?

\providecommand{\bysame}{\leavevmode\hbox to3em{\hrulefill}\thinspace}
\providecommand{\MR}{\relax\ifhmode\unskip\space\fi MR }
\providecommand{\MRhref}[2]{%
  \href{http://www.ams.org/mathscinet-getitem?mr=#1}{#2}
}
\providecommand{\href}[2]{#2}

\end{document}